\documentclass[11pt]{amsart}
\usepackage{epsfig}
\usepackage{amsthm}
\usepackage{amssymb}
\usepackage{amsmath}
\usepackage{graphicx}

\newtheorem {Theorem}   {Theorem}

\numberwithin{Theorem}{section}

\theoremstyle{definition}

\theoremstyle{remark}
\newtheorem{Remark}[Theorem]{Remark}

\expandafter\chardef\csname pre amssym.def
at\endcsname=\the\catcode`\@ \catcode`\@=11
\def\undefine#1{\let#1\undefined}
\def\newsymbol#1#2#3#4#5{\let\next@\relax
 \ifnum#2=\@ne\let\next@\msafam@\else
 \ifnum#2=\tw@\let\next@\msbfam@\fi\fi
 \mathchardef#1="#3\next@#4#5}
\def\mathhexbox@#1#2#3{\relax
 \ifmmode\mathpalette{}{\m@th\mathchar"#1#2#3}%
 \else\leavevmode\hbox{$\m@th\mathchar"#1#2#3$}\fi}
\def\hexnumber@#1{\ifcase#1 0\or 1\or 2\or 3\or 4\or 5\or 6\or 7\or 8\or
 9\or A\or B\or C\or D\or E\or F\fi}

\font\teneufm=eufm10 \font\seveneufm=eufm7 \font\fiveeufm=eufm5
\newfam\eufmfam
\textfont\eufmfam=\teneufm \scriptfont\eufmfam=\seveneufm
\scriptscriptfont\eufmfam=\fiveeufm

\catcode`\@=\csname pre amssym.def at\endcsname

\def    \C {{\mathbb C}}
\def    \reals  {{\mathbb R}}
\def    \R      {{\mathbb R}}

\def    \p      {\partial}

\begin{document}





\title[The N-Vortex Problem on
a Symmetric Ellipsoid ]{The N-Vortex Problem on a Symmetric
Ellipsoid: A Perturbation Approach.}

\author[Castilho and Machado]{C\'esar Castilho$^{\dagger\, \star}$ and H\'elio Machado$^{\dagger}$}
\maketitle
 \centerline {$^{\dagger}$Departamento de Matem\'atica} \par
\centerline{ Universidade Federal de Pernambuco} \par
\centerline{Recife, PE CEP 50740-540 Brazil} \vspace{0.2cm}
 \centerline{and} \vspace{0.2cm}
\par \centerline{$^{\star}$ The
Abdus Salam ICTP} \par \centerline{Strada Costiera 11 } \par
\centerline{ Trieste 34100 Italy}


\bigskip

\begin{abstract}
We consider the N-vortex problem on a  ellipsoid of revolution.
Applying standard techniques of classical perturbation theory we
construct a  sequence of conformal transformations from the
ellipsoid into the complex plane.  Using these transformations the
equations of motion for the N-vortex problem on the ellipsoid are
written as a formal series on the eccentricity of the ellipsoid's
generating ellipse. First order equations are obtained explicitly.
We show numerically that the truncated first order system for the
three-vortices system on the symmetric ellipsoid is
non-integrable.

\end{abstract}
\vskip 1.0cm

\centerline{ {\bf Key Words:} Vortex Dynamics, Hamiltonian
Systems, Perturbation Theory.}\par \vskip 1.0cm \centerline{ {\bf
MSC:} 76C05, 58F40, 70H05.}

\section{Introduction}
The equations of point vortices motion in the plane were
introduced by Helmholtz  \cite{helm} and described as an
Hamiltonian system by Kirchhoff in 1876 \cite{kirch}. These
equations were first generalized to describe the motion of
point-vortices on a sphere by Bogolmonov in 1977 \cite{bol}.
Bogomolnov's work was put into solid mathematical grounds by
Kimura and Okamoto in 1987 \cite{kim}.  Independently Hally
 wrote in 1979 \cite{hal} the equations for vortices motion in
symmetric surfaces of revolution (which englobe the sphere). To
write his equations Hally represented the surface $M$ on the
complex plane with a metric conformal to the Euclidean metric: Let
the $i$-th vortex, with vorticity $\Gamma_i$ be represented by the
complex coordinate $\xi_i$. Calling the surface's conformal factor
$h (\xi \, , \bar \xi)$, Hally's equations are
\begin{equation}
\label{hally0} {\bar \xi_n}^{\, \prime} = h^{-2} (\xi_n, \bar
\xi_n) \left(\sum_{k}^{\prime} -i \frac{\Gamma_k}{\xi_n - \xi_k} +
i \, \Gamma_n \, \frac{\p}{\p \xi_n}\log ( \, h(\xi_n,\bar \xi_n)
\, ) \right) \, .\end{equation} The difficulty in working with
these equations is that one must know explicitly the conformal
factor $h (\xi \, , \bar \xi)$ and this factor is not known for
many surfaces but the sphere.  \par

The aim of this paper is to write the equations for the N-vortex
problem on the ellipsoid of revolution
$$
\label{ellipsoid}
  \frac{{x}^2}{R^2} + \frac{{y}^2}{R^2} + \frac{{z}^2}{R^2\, (1+\epsilon)}
   = 1
$$
as an $\epsilon$-series when $|\epsilon| << 1 $. To achieve this
goal we write a perturbative  coordinate transformation from the
ellipsoid into the plane (see (\ref{trans}))and impose
conformality at each order. It turns out that the transformation's
n-th order term is given by the solution of a certain linear
differential equation and that the first order term can be
computed explicitly (see (\ref{equapert})). Using Hally's
equations we calculated rigorously the first order perturbation
term  for the vortex equations on the ellipsoid. In spherical
coordinates the truncated first order equations are given by the
hamiltonian system with hamiltonian
\begin{equation} \label{central}
 H = H_0 + \epsilon \, H_1
\end{equation}
where $$H_0 = \frac{1}{2} \, \sum_{i,j}^{\prime} \Gamma_i \,
\Gamma_j \, \left[ \ln \left( 2 - 2 \, \cos(\theta_i) \,
\cos(\theta_j) - 2 \, \sin(\theta_i) \, \sin(\theta_j) \,
\cos(\phi_i - \phi_j) \right) \right] \, ,$$
$$ H_1 =   \frac{1}{2} \, \sum_{i,j}^{\prime} \Gamma_i \,
\Gamma_j \, \left[  \left( \cos(\theta_i)^2
 + \cos(\theta_j)^2 \right) \right]$$
 and symplectic two-form
$$ w = \sum_{i=1}^N \Gamma_i \, \sin(\theta_i) \, d \theta_i \wedge d \phi_i  \, . $$

The paper is organized as follows: In section (\ref{hallysec}) we
introduce all the concepts used in the derivation of the vortex
equations. Hally's equations are written as an invariant
hamiltonian system and, for the spherical case, are shown to be
equivalent to Bogolmonov's equations. In section
(\ref{pertubative}) we derive the expression for the first order
perturbation term. The main idea is to write the conformal
transformation from the ellipsoid into the plane  as a formal
series in the eccentricity of the generating ellipse of the
ellipsoid. In section (\ref{aplicacao}) we exhibit two simple
applications of  equations (\ref{central}). The case $N=2$,
despite being integrable, seems not to admit closed solutions. We
show that, when $|\Gamma_1|=|\Gamma_2|$, the only solutions for
which the two vortices maintain their relative distance constant
are the relative equilibria. For the $N=3$ case, we show
numerically, by choosing specific parameters, that the
three-vortices problem (\ref{central}) on the ellipsoid is
chaotic.
\par

\section{The Equations of Hally}
\label{hallysec}
 In this section we discuss equations (\ref{hally0}) from a geometric
and invariant point of view. Hally's equations were deduced under
the topological constraint that the sum of all  vorticities is
zero and, in fact, this is essential in recasting the equations in
Hamiltonian form. We also show, that, under this zero vorticity
constraint,  Hally's equations on the sphere, coincide with
Bogolmonov's equations.
\subsection{Hally's equations as an hamiltonian system on the conformal plane} Let $M$ be a
2-dimensional Riemannian manifold with metric tensor $g$. $M$ is a
symplectic manifold with symplectic form $w$ given by the area
form induced by $g$. Let $H : M \rightarrow \reals $ be the
hamiltonian function. Hamilton's equations \cite{marsden } are
\begin{equation}
\label{hamil}
 w ( X_H \, , \, \cdot) = dH \,
 \end{equation}
where $X_H$ is the hamiltonian vector field induced by $H$. Now we
assume that $M$ is conformal to the plane, that is, that there are
coordinates $x_i \colon M \to \reals^2 \, ; \, i=1,2 $ where the
metric is given by
$$ g \, (x_1,x_2) = h^2(x_1,x_2) \,\left( dx_1 \otimes dx_1 + dx_2 \otimes dx_2 \right) $$
for  some function $h$ . In these coordinates the symplectic (area)
form is given by
$$ w = h^2(x_1 \, , \, x_2) \, dx_1 \wedge dx_2 \, . $$
Introduce complex coordinates $(\xi \, , \, \bar \xi)$ on the
plane $(x_1 \, , \, x_2)$  by
$$ \begin{array}{l}
\xi= \frac{1}{\sqrt{2}}\left( x_1 + i \, x_2 \right) \, , \\
\\
\bar \xi = \frac{1}{\sqrt{2}}\left( x_1 - i \, x_2 \right) \, .
\end{array} $$
Using $(\xi \, , \, \bar \xi)$ we obtain
$$g(\xi,\bar \xi) =   h^2(\xi \, , \, \bar \xi) \, \left( d\xi \otimes d \xi + d{\bar \xi} \otimes d{\bar \xi} \right)
\,  $$ and
$$ w = - i \, h^2(\xi, \bar \xi) \, d \xi \wedge d {\bar \xi} \, . $$
Hamilton's equations (\ref{hamil}) become
$$    h^2(\xi, \bar \xi) \, d\xi \wedge d {\bar \xi}  \; ( X_H \, , \, \cdot)= d \left( i \, H \right)  \, .$$
For convenience we think of $ h^2(\xi, \bar \xi) \, d\xi \wedge d
{\bar \xi} $ as the new area form and $ i \,  H $ as the new
hamiltonian function.
 Writing $X_H = \dot \xi \,  \frac{\p}{\p
\xi} + \dot {\bar \xi} \frac{\p}{\p \bar \xi}$ the equations
become
\begin{equation}
\label{equa1} h^2(\xi, \bar \xi) \, \dot {\bar \xi} =- i \,
\frac{\p H}{\p \xi} \, .  \end{equation}

We now show that Hally's equations can be written as a hamiltonian
system on the conformal plane. The equations are
\begin{equation}
\label{hally1} {\bar \xi_n}^{\, \prime} = h^{-2} (\xi_n, \bar \xi_n)
\left(\sum_{k}^{\prime} -i \frac{\Gamma_k}{\xi_n - \xi_k} + i \,
\Gamma_n \, \frac{\p}{\p \xi_n}\log ( \, h(\xi_n,\bar \xi_n) \, )
\right) \, .\end{equation}
 Multiplying  each term of the sum in
the right hand side by $\displaystyle \frac{\bar \xi_n - \bar
\xi_k }{{\bar \xi_n - \bar \xi_k}^{\phantom{*}} }$ and adding the
zero $ \displaystyle i \, \Gamma_n \, \frac{\p}{\p \xi_n}\, \log
(h(\xi_k \, , \, \bar \xi_k))\, $ to the right-hand side we obtain
$$ \begin{array}{lll} {\bar \xi_n}^{\, \prime} &  = & h^{-2} (\xi_n, \bar \xi_n) \left({\sum_k}^{\prime} -i \,
 \Gamma_k \frac{\p}{\p \xi_n}{\log |\xi_n - \xi_k|^2} + \right. \\
 \\
 & &
\hskip 6.0cm \left. i \, \Gamma_n \, \frac{\p}{\p \xi_n} \log ( \,
h(\xi_n,\bar \xi_n) \, h(\xi_k,\bar \xi_k) \, )  \ \right) \, .
\end{array} $$ Assuming that the sum of the vorticities is equal to zero we
have
$$ \begin{array}{lll} h^2 (\xi_n, \bar \xi_n)\, {\bar \xi_n}^{\, \prime}&  = & -i \,  \left({\sum_\xi}^{\prime}\,
 \Gamma_k \, \frac{\p}{\p \xi_n}{\log  |\xi_n - \xi_k |^2} + \right. \\
 \\
 & & \hskip 4.0cm \left. \sum_k^{\prime} \, \Gamma_k \, \, \frac{\p}{\p \xi_n}\log ( \, h(\xi_n,\bar \xi_n) \, h(\xi_k,\bar \xi_k)\, ) \right)
\, , \end{array}$$ that is

\begin{equation}
\label{hally2}
 h^2 (\xi_n, \bar \xi_n)\, {\bar \xi_n}^{\, \prime}  =  -i \,
{\sum_k}^{\prime}\,
 \Gamma_k \, \frac{\p}{\p \xi_n}{\log \left( h(\xi_n,\bar \xi_n) \,  h(\xi_k,\bar \xi_k) |\xi_n - \xi_k|^2
 \right)}\, .
\end{equation}
But those are Hamilton's equations
\begin{equation}
\label{hamilton} w (X_H, \cdot) = d (i H)  \end{equation} for the
symplectic two-form
\begin{equation}
\label{forma} \tilde w = \sum_{i=1}^N \Gamma_i \, h^2(\xi_i \, ,
\, \bar \xi_i) \, d \xi_i \wedge d \bar \xi_i \,  \end{equation}
and Hamiltonian function
\begin{equation}
\label{hamilconfor}
 H = \frac{1}{2}\,{\sum_{k\, ,\, n}}^{\prime}\,
 \Gamma_k \, \Gamma_n {\log \left( h(\xi_n,\bar \xi_n)
  \,  h(\xi_k,\bar \xi_k) |\xi_n - \xi_k |^2 \right)} \, .\end{equation}

\subsection{From Hally's to Bogolmonov's equations}
 Denote by $S^2
\subset \reals^3 $ the two-dimensional sphere
$$ S^2 = \left\{(x,y,z) \in \reals^3 \, | \, x^2 + y^2 + z^2 =
R^2 \right\}\, . $$ $S^2$ is a Riemannian manifold with the metric
induced by the Euclidean metric in $\reals^3$. Introduce spherical
coordinates on $S^2$ trough the relations
$$\begin{array}{l}  x = R \cos(\phi) \, \sin(\theta) \, ,\\
                    y = R \sin(\phi) \, \sin(\theta) \, ,\\
                    z = R \, \cos(\theta) \, , \end{array} $$
                    with $ 0
\le \phi < 2 \, \pi $, $0 <  \theta < \pi \, .$

Consider the conformal map  $ \sigma : S^2 \rightarrow \C $,
$\sigma(\theta \, , \, \phi) = u + i \, v $ where
\begin{equation}
\label{conf}
\begin{array}{l}
 u =\tan(\theta /2) \, \cos(\phi) \, , \\
 \\
 v = \tan(\theta /2) \, \sin(\phi)  \, .
 \end{array}
 \end{equation}
The conformal factor is given by
$$ h(\xi,\bar \xi) = \frac{2 \, R }{1 + \xi \, \bar \xi} \,  $$
where $\xi = u + i \, v \, . $ Now let $\vec r_1 \, , \,  \vec r_2
\in S^2 \subset \reals^3 $. The Euclidean distance $\| \vec r_1 -
\vec r_2 \|^2 $ in $\reals^3$ can be computed in the conformal
variables. In fact
$$\begin{array}{rll} \| \vec r_1 -
\vec r_2 \|^2 & = &  (x_1 - x_2)^2 + (y_1 - y_2)^2 + (z_1 -z_2)^2
\,
,\\
\\
& = &  R^2 \, \left\{2 - 2 \, \cos(\theta_1) \, \cos(\theta_2) -
2\, \sin(\theta_1) \, \sin(\theta_2)\,
\cos(\phi_1 - \phi_2)  \right\} \, ,\\
\\
 & = & 4 \, \cos(\theta_1)^2 \, \cos(\theta_2)^2 \, R^2
\,\left\{(u_1^2 + v_1^2) + (u_1^2 + v_1^2) - 2 \, u_1 \, v_1 - 2
\, u_2 \, v_2 \right\} \, , \\
\\
& = & \frac{2 \, R}{1 + \xi_1 \, \bar \xi_1} \, \frac{2 \, R}{1 +
\xi_2 \, {\bar \xi_2}}^{\phantom{*}} \, | \xi_1 - \xi_2 |^2 \, ;
\end{array}
$$
that is
$$\| \vec r_1 -
\vec r_2 \|^2 =
 h(\xi_1 ,\bar \xi_1) \,  h(\xi_2 ,\bar \xi_2)\, \, | \xi_1 - \xi_2 |^2 \, .$$

Now consider the map $F: S^2 \times \dots \times S^2 \rightarrow
\C^N \, $ given by
$$ F(\vec r_1 \, , \, \vec r_2 \, , \dots \, , \, \vec r_N )
= ( \sigma (\vec r_1) \, , \, \sigma(\vec r_2) \, , \dots \, , \,
\sigma(\vec r_N) ) \, . $$ The pull back by $F$ of the symplectic
form $$ w = \sum_{i=1}^N \Gamma_i \, h^2(\xi_i \, , \, \bar \xi_i)
\, d \xi_i \wedge d \bar \xi_i
$$
is just
$$ F^{*} \tilde w = \sum_{i=1}^N \Gamma_i \, w_i \, $$
where by $w_i$ we denote the standard area form
 \begin{equation}
 \label{simp} w_i =  x_i \, dy_i \wedge dz_i + y_i \, dz_i \wedge dx_i + z_i
 \, dx_i \wedge dy_i \,
 \end{equation}
induced by the Euclidean metric on $S_i^2$. The pull back of the
Hamiltonian function

$$ H = \frac{1}{2}\,{\sum_{k\, ,\, n}}^{\prime}\,
 \Gamma_k \, \Gamma_n {\log \left( h(\xi_n,\bar \xi_n)
  \,  h(\xi_k,\bar \xi_k) \|\xi_n - \xi_k\|^2 \right)} $$
is

 \begin{equation}
 \label{hamm} F^{*} H = \frac{1}{2}\,{\sum_{k\, ,\, n}}^{\prime}\,
 \Gamma_k \, \Gamma_n {\log \left( \| \vec r_k -
\vec r_n \|^2 \right)} \, .\end{equation} Hamilton's equations for
Hamiltonian function (\ref{hamm}) and symplectic form (\ref{simp})
are precisely Bogolmonov's equations for $N$ interacting vortices
on a sphere.

We notice that Hally's equations are equivalent to Bogolmonov's
equations only for the case where all the vorticities sums up to
zero. Therefore, in dealing with the general vorticity case, one
must use Bogomolnov's equations.

\section{Equations on the Symmetric Ellipsoid - Perturbation Approach}
\label{pertubative}
 Let   $E^2 \subset \R^3$ represent
 the two dimensional ellipsoid of revolution
\begin{equation}
\label{ellipsoid}
  \frac{{x}^2}{R^2} + \frac{{y}^2}{R^2} + \frac{{z}^2}{R^2\, (1+\epsilon)}
   = 1 \, .
\end{equation}
Where $ R > 0 $ and $|\epsilon | << 1 \, . $ If $\epsilon  > 0 $
the ellipsoid is prolate. In this case the eccentricity $e$ of the
generating ellipse is given by $e= \sqrt{\frac{\epsilon}{(1 +
\epsilon)}}\,  $. If $\epsilon < 0 $ the ellipsoid is oblate and
the eccentricity of the generating ellipse is given by $ e =
\sqrt{-\epsilon}\, \, . $ We observe that  $\epsilon= \mathcal O
(e^2)$ and therefore the truncated first order system is expected
to approximate the real system for quite high values of the
eccentricity. \par

 We introduce coordinates on $E^2$ trough
the relations
\begin{equation}
\label{coord11}
\begin{array}{l}  x = R \cos(\phi) \, \sin(\theta) \, ,\\
                    y = R \sin(\phi) \, \sin(\theta) \, ,\\
                    z = R \, \sqrt{ 1 + \epsilon} \, \cos(\theta) \, , \end{array}
                    \end{equation}
                    with $ 0
\le \phi < 2 \, \pi $, $0 <  \theta < \pi \, .$
 If $\epsilon = 0 $
the surface is a sphere of radius $R$. In this case the
 projection
\begin{equation}
\begin{array}{l}
\label{stereo}
 u = \tan(\theta /2) \, \cos(\phi) \, , \\
 \\
 v = \tan(\theta /2) \, \sin(\phi) \,
\end{array}
\end{equation}
is a conformal transformation. We want to find a conformal
transformation from $E^2$ to the conformal plane when $\epsilon
\ne 0 $. By the symmetry of the ellipsoid the conformal factor
must depend only on the radial variable $r$ and not on the angular
variable $\psi \, . $ We look for a transformation of the form
\begin{equation}
\begin{array}{l}
\label{trans}
 \theta = 2 \, \arctan(r) + 2 \, \sum_{i=1}^{\infty} f_i(r) \, \epsilon^i \, , \\
 \\
 \phi = \psi \, .
\end{array}
\end{equation}

We observe that when $\epsilon = 0 $ the inverse of the above
transformation is just transformation (\ref{stereo}) written in
polar coordinates. \par

 The symmetry of the ellipsoid implies that the
equator is a circle of radius $R$. In the $(\theta, \phi)$
coordinates the equator has parametrization given by  $\theta =
\frac{\pi}{2}$ and $\phi \in [ 0 , 2 \, \pi )\, $. Therefore in
the conformal coordinates $(u,v)$ the equator becomes the circle
of radius $ r = u^2 + v^2 = 1 \, . $ We remark that the conformal
factor for the ellipsoid, when restricted to the circle of radius
$1$  is just the usual one for the sphere since the ellipsoid
$E^2$ and the $S^2$ coincide at the equator. This is of utmost
importance for what follows. In fact, since we will impose
conformality at each order this trivial observation implies that
$f_i(1) = 0 $ for all $i$. This is the condition that will
guarantee the uniqueness of our series expansion.
\begin{Remark}
\label{unique} Alternatively (and equivalently), the uniqueness of
the expansion can be obtained by imposing that the symmetry
condition $ h(r) = \frac{1}{r} \, h(r^{-1}) $ has to be satisfied
at all orders.
\end{Remark}

 The algorithmic determination of the $f_i(r)$
 is obtained by imposing conformality
at each order  and follows the standard procedures of classical
perturbation theory: Truncate the series at order $n$ , find the
conditions to be satisfied, and them assume those conditions to
find the term of order $n+1$. To find $f_1(r)$ we impose that
$$ \begin{array}{l}
 \theta = 2 \arctan(r) + 2 \, f_1(r)\, \epsilon \, , \\
 \\
 \phi = \psi \, .
\end{array} $$
generates a conformal transformation up to first order in
$\epsilon$. \par The line element for the ellipsoid is
 $$ ds^2
= R^2 \sin^2 (\theta) \, d \phi^2 + R^2 \, (1 +
  \epsilon \, \sin^2(\theta))\, d\theta^2 \,
 . $$
We want to write $ds^2$ in terms of $r$ and $\phi$. To this goal
it suffices writing $\sin(\theta)$ and $d\theta$ in terms of $r$.
But
$$\begin{array}{l}
\sin(\theta) =\sin  \left( 2 \, \arctan(r) + 2 \, f_1(r)\,
\epsilon
\right) \, , \\
\\
= \sin  ( 2 \, \arctan(r)) + 2 \, \cos( 2 \, \arctan(r)) \, f_1(r)
\, \epsilon + \mathcal
O(\epsilon^2) \, \\
\\
= 2 \,\sin  (\arctan(r))\, \cos(\arctan(r)) + \\
\\
+ 2 \left( \cos( \arctan(r))^2 - \sin  (\arctan(r))^2 \right) \,
f_1(r)  \, \epsilon + \mathcal O(\epsilon^2) \, , \\
\\
= 2 \, \frac{r}{1+r^2} + 2 \frac{(1-r^2)}{1+r^2} \, f_1(r) \,
\epsilon +  \mathcal O(\epsilon^2) \, .
\end{array}
$$
This gives
$$ \sin(\theta)^2 = \frac{4 \, r^2}{(1+r^2)^2} + \frac{8 \, r\, (1-r^2)}{(1+r^2)^2}\, f_1(r) \,
\epsilon +  \mathcal O(\epsilon^2) \, .$$ Also
 $$ d \theta = \left( 2 \frac{1}{1+r^2} + 2 \, \frac{\p f_1(r)}{\p
 r} \, \epsilon +  \mathcal O(\epsilon^2) \right) \, \, dr \, .$$
 Therefore
 $$d \theta^2 = \left( 4 \frac{1}{(1+r^2)^2}  + \frac{8}{1+r^2} \, \frac{\p f_1(r)}{\p
 r}\, \epsilon \, +  \mathcal O(\epsilon^2)  \right) \, dr^2 \, .$$
Finally
$$
\begin{array}{r} ds^2=  \frac{4 \, R^2}{(1+r^2)^2}\left\{ \left[ 1+
\frac{2 \,(1-r^2)}{r}\, f_1(r) \, \epsilon  \, +  \mathcal
O(\epsilon^2) \right] \, r^2 d \psi^2  + \right.
\\
\\
+ \left. \left[ 1 + \left( \frac{4 \, r^2}{(1+r^2)^2} + 2 \, (1 +
r^2) \, \frac{\p f_1(r)}{\p r} \right) \, \epsilon  \, +  \mathcal
O(\epsilon^2) \right] \, dr^2 \,\right\} .
\end{array}
$$
This will be conformal to the plane if it is a multiple of $ r^2
\, d\psi^2 + dr^2 $. We observe that if $\epsilon = 0 $ then the
transformation is conformal with conformal factor $ \sqrt{\frac{4
\, R^2}{(1+r^2)^2}}$. For $\epsilon \ne 0 $ the only way this
transformation can be conformal for all values of $\epsilon$ is if
the coefficient of $\epsilon$ in the first square bracket is equal
to the coefficient of $\epsilon$ in the second square bracket,
that is if
\begin{equation}
\label{equapert}
 (1 + r^2) \, \frac{\p
f_1(r)}{\p r} - \frac{(1-r^2)}{r}\, f_1(r) = -\frac{2 \,
r^2}{(1+r^2)^2} \, .$$ This is a linear differential equation for
$f_1(r)$. Its general solution (for $r > 0 $) is equal to
$$ f_1(r) = -\frac{( \frac{1}{1+r^2}+ c )}{(1+r^2)} \, r \, .
\end{equation}
 Now by remark (\ref{unique}) we must have $f_1(1)=0$ what implies
 $c=-\frac{1}{2}$. Finally we have
 $$ f_1(r) =  \frac{1 - r^2}{(1+r^2)^2}\, r \, . $$
Transformation (\ref{trans}) becomes
\begin{equation}
\begin{array}{l}
\label{trans1}
 \theta = 2 \, \arctan(r) + \frac{(1-r^2)}{(1+r^2)^2}\, r \, \epsilon+  \sum_{i=2}^{\infty} f_i(r) \, \epsilon^i \, , \\
 \\
 \phi = \psi \, .
\end{array}
\end{equation}
The conformal factor is, to first order given by,
$$ h(r)^2 = \frac{4 \, R^2}{(1+r^2)^2}\, \left\{ 1 + 2 \, \frac{(1-r^2)^2}{(1+r^2)^2} \, \epsilon \right\}  $$
that is
$$ h(r) =\frac{2 \, R}{(1+r^2)}\, \left\{ 1 +\frac{(1-r^2)^2}{(1+r^2)^2}\, \epsilon
\right\}  \, $$ We write
\begin{equation}
\label{expan}
 h(r) = h_0(r) + h_1(r) \,
\epsilon + \mathcal O(\epsilon^2) \, ,\end{equation}
 where
$$h_0(r) =\frac{2 \, R}{(1+r^2)} \, \, \, \mbox{and}
\, \, \, h_1(r) = 2 \, R \,  \frac{(1-r^2)^2}{(1+r^2)^3} \, .$$

For future reference we observe that in terms of the variables
$\theta$ and $\phi$ we obtain
$$h_0(r) =2 \, R \cos^2( \theta /2) \, \, \, \mbox{and}
\, \, \, h_1(r) = 4 \, R \, \cos(\theta) \, \cos^2(\theta /2) \, .$$

\begin{Remark} The above expansion can not be carried
explicitly further. In fact, the equation for the second order
factor will be
$$ {\frac {d \, f_2}{dr}} + \frac{1}{2} {\frac {r
\left( {r}^{2} - 1  \right)  }{ r^2+1}} \, f_2   +{\frac {1-13\,{r
}^{2}+46\,{r}^{4}-{r}^{10}-26\,{r}^{6}+9\,{r}^{8}}{4 \,  \left(
1+{r}^{2}
 \right) ^{5}}} = 0 \, .
 $$
The general solution for this equation is given by
$$
f_2(r) =\frac{\sqrt {1+{r}^{2}}}{4 \, {e^{{r}^{2}/4}}}\, \left( \int
\!{\frac { \left(( {r}^{2} - 1)^2 + 4 \, r^2 \right) \,  \left(
{r}^{6}-3\,{r}^ {4}+7\,{r}^{2}-1 \right) }{ \left( 1+{r}^{2}
 \right) ^{11/2}}} \, {e^{1/4\,{r}^{2}}} \, {dr}+ 4\,c \right)
$$
where $c$ is an arbitrary constant. The above integral can not be
computed using elementary functions.
\end{Remark}

\subsection{N-Vortex Equations:}Using the conformal factor just found
we will write the N-vortex equations on the ellipsoid of
revolution up to first order terms in $\epsilon$. We use the
coordinates $(\theta, \phi)$ introduced trough equations
(\ref{coord11}). The symplectic two-form becomes
$$ w = \sqrt{1 + \epsilon} \, \; d \phi \wedge d (\cos(\theta)) \, .$$
Observing that
$$  \sqrt{1 + \epsilon}  \; d \phi \wedge d (\cos(\theta)) \left( X_H \, , \, \cdot \, \right)\, =
 d \phi \wedge d (\cos(\theta)) \left( \sqrt{1 + \epsilon} \, \, X_H \, , \, \cdot \, \right)\, ,$$
we define the vector field $\tilde {X} = \sqrt{1 + \epsilon} \, \,
X_H $.  $\tilde {X}$ satisfies Hamilton's equations
 \begin{equation}
 \label{ham}
  d \phi \wedge d (\cos(\theta)) \left( \tilde{X} \, , \, \cdot \,
 \right)\, = dH \, , \end{equation}
and the flow  $\phi_s(x)$ of the vector field $\tilde{X_h}$ is a
time reparametrization of the flow $\phi_t(x)$induced by $ X_H $. In
fact $ t= \sqrt{1 + \epsilon} \; s . $ We therefore consider
equation (\ref{ham}) with $H$ and $h$ given by (\ref{hamilconfor})
and  (\ref{expan}) respectively . $H$ is given by
$$ H =\frac{1}{2} \, {\sum_{i,j}}^{\prime} \Gamma_i \, \Gamma_j \, \ln
\left[ \left(h_0(r_i) + \epsilon \, h_1(r_i) \right) \,
\left(h_0(r_j) + \epsilon \, h_1(r_j) \,  \right)\, | \xi_i - \xi_j
|^2 + \mathcal O (\epsilon^2) \right]
$$
that is
$$H =\frac{1}{2} \, {\sum_{i,j}}^{\prime} \Gamma_i \, \Gamma_j \, \left[ \ln
\left( h_0(r_i) \, h_0(r_j)\, | \xi_i - \xi_j |^2  \right) +
\epsilon \, \left( \frac{h_1(r_i)}{h_0(r_i)} +
\frac{h_1(r_j)}{h_0(r_j)}\right) \right] + \mathcal O (\epsilon^2)
 $$
 The Hamiltonian function $H$ is a
function on the variables $(\theta,\phi)$. Moreover, the
symplectic two-form (\ref{ham}) is the standard symplectic form
over the sphere. Therefore, we write

$$\begin{array}{r} H =\frac{1}{2} \, {\sum_{i,j}}^{\prime} \Gamma_i \, \Gamma_j \, \left[ \ln
\left( 2 - 2 \, \cos(\theta_i) \, \cos(\theta_j) - 2 \,
\sin(\theta_i) \, \sin(\theta_j) \, \cos(\phi_i - \phi_j) \right)
\right.
\\
\\
\left. + \epsilon \, \left( \cos(\theta_i)^2
 + \cos(\theta_j)^2  \right) \right] + \mathcal O
 (\epsilon^2)\,.
\end{array} $$
We conclude that, to first order in $\epsilon$ the dynamics of the N
vortices on the ellipsoid  is determined by the hamiltonian system
with Hamiltonian function
$$\begin{array}{r} H =\frac{1}{2} \, \sum_{i,j}^{\prime} \Gamma_i \, \Gamma_j \, \left[ \ln
\left( 2 - 2 \, \cos(\theta_i) \, \cos(\theta_j) - 2 \,
\sin(\theta_i) \, \sin(\theta_j) \, \cos(\phi_i - \phi_j) \right)
\right.
\\
\\
\left. + \epsilon \, \left( \cos(\theta_i)^2
 + \cos(\theta_j)^2 \right) \right]
\end{array} $$
and symplectic two-form
$$ w = \sum_{i=1}^N \Gamma_i \, d \phi_i \wedge d (\cos(\theta_i)) \, . $$
Alternatively, in cartesian coordinates the above Hamiltonian
system becomes
\begin{equation}
\label{hamN}
 H =
\frac{1}{2} \, \sum_{i,j}^{\prime} \Gamma_i \, \Gamma_j \, \log \|
\vec r_i - \vec r_j \|^2 + \frac{\epsilon}{2} \, \sum_{i,j}^{\prime}
\Gamma_i \, \Gamma_j (z_i^2 + z_j^2) \, ,
\end{equation}
with symplectic two-form given by
\begin{equation}
\label{simpN}
 w = \sum_{i=1}^n
\Gamma_i \, \left( x_i \, dy_i \wedge dz_i
 +  y_i \, dz_i \wedge dx_i  +  z_i \, dx_i \wedge dy_i\, \right) \, .
 \end{equation}
 \section{Applications: The N=2,3 cases.}
\label{aplicacao}
\subsection{Two-Vortex problem}
When $N=2$ the hamiltonian system given by (\ref{hamN}) and
(\ref{simpN}) reduces to
\begin{equation} \label{mot1}
\begin{array}{l}
         \dot{\vec r }_1 \, = \, \Gamma_2 \,
         \frac{\vec r_2 \, \times \, \vec r_1}{\| \vec r_1 - \vec r_2
         \|^{2}}  + \epsilon \, \Gamma_2 \, \vec r_1 \times z_1 \, \hat z \, , \\
         \phantom{x} \\
         \dot{\vec r }_2 \, = \, \Gamma_1 \,
         \frac{\vec r_1 \, \times \, \vec r_2}{\| \vec r_1 - \vec r_2
         \|^2} + \epsilon \, \Gamma_1 \, \vec r_2 \times z_2\, \hat z
         \, .
          \end{array}
\end{equation}
Doing the following change of variables
\begin{equation}
\label{cha1}  \begin{array}{l}
 \vec w = \Gamma_1 \, \vec {r_1} - \Gamma_2 \, \vec
{r _2} \, , \\
\phantom{x} \\
\vec c = \Gamma_1 \, \vec {r_1} + \Gamma_2 \, \vec {r _2} \,
.\end{array}  \end{equation} we obtain
\begin{equation}
\label{mot2} \begin{array}{l}
 \dot {\vec w} = \frac{\vec c \, \times \, \vec w}{\| \, \alpha_1 \, \vec
c \, - \, \alpha_2 \, \vec w \, \|^{2}}\, +  \epsilon \,   \left(
\Delta_1(c_z,w_z) \, \vec c \times \hat z
+ \Delta_2(c_z,w_z) \, \vec w \times \, \hat z \right)\\
\phantom{x} \\
\dot{ \vec c} = \epsilon \, \left( \Delta_2(c_z,w_z) \, \vec c
\times \hat z + \Delta_1(c_z,w_z) \, \vec w \times \, \hat z
\right) \, ;\end{array}
\end{equation}
where
$$ \alpha_1 \equiv \frac{1}{2} \, \frac{\Gamma_2 - \Gamma_1}{\Gamma_1 \,
\Gamma_2}  \;\;\;\;\; , \;\;\;\;\;\alpha_2 \equiv \frac{1}{2} \,
\frac{\Gamma_1 + \Gamma_2}{\Gamma_1 \, \Gamma_2}\, ,
$$
$$ \Delta_1 \equiv \frac{c_z + w_z}{4 \, \Gamma_1^2} + \frac{c_z - w_z}{4 \, \Gamma_2^2}
 \;\;\;\;\; \mbox{and} \;\;\;\;\; \Delta_2 \equiv \frac{c_z + w_z}{4 \, \Gamma_1^2} - \frac{c_z - w_z}{4 \, \Gamma_2^2}.
$$ Equations (\ref{mot2}) are singular when $\alpha_1 \, \vec c \, =
\, \alpha_2 \, \vec w \, ,$ i.e. $\vec r_1 = \vec r_2 \,$.
 When $\epsilon = 0 $ equations (\ref{mot2}) reduces to
a linear system and can be easily solved. Those are the well known
solutions for the two-vortex problem over the sphere: In a system
of coordinates such that its z-axis is in the same direction as
the center of vorticity $\vec c$ the solutions are given in the
$\vec r_1$, $\vec r_2$ variables by
\begin{equation}
\label{2vort}
\begin{array}{l}
\vec{r_1}(t) =\frac{1}{2 \, \Gamma_1}\left( \mu  \, \cos \left(
\theta \right) \, , \, \mu  \, \sin \left( \theta \right) \, , \, \| \,\vec c \, \|+ w_z \right)\, , \\
\phantom{x} \\
 \vec{r_2}(t) =\frac{1}{2 \, \Gamma_2}\left( -\mu  \, \cos \left(
\theta \right) \, , \, -\mu  \, \sin \left( \theta \right) \, , \,
\| \,\vec c \, \| - w_z \right)\, \end{array}
\end{equation}
where
$$ \theta \equiv \gamma \, t + \phi \, , $$
with $\phi$ a constant and
\begin{equation}
\label{gamma}
 \gamma
= \frac{ \| \,\vec c \, \|}{\alpha_2^2\,\mu ^2 + \left(\alpha_2\,
w_z - \alpha_1 \,  \| \,\vec c \, \| \right)^2}\, ,\end{equation}
with
$$ \mu^2 =w_x^2 + w_y^2 = w^2 -w_z^2 \,  \, \, \mbox{and} \, \, \,  \mu \ge 0 .$$
We observe that all solutions are periodic and the relative
distance between the vortices is a first integral.
\par
 When $\epsilon \ne 0$  the center of vorticity vector $\vec c$ is not preserved
 . In fact, only it's $z$ component $c_z$ is
preserved. This first integral implies, by dimensionality that
system is integrable. Since $\| \vec r_1 \| = \| \vec r_2 \| = 1 $
we obtain that
$$ \| \vec c \|^2 + \| \vec w \|^2 = k_1 \, \, \, \, \mbox{and} \, \, \, \,
\vec c \cdot \vec w = k_2 \,
$$ are constants of motion. Despite being integrable system
(\ref{mot2}) seems not to admit explicit solutions for a general
initial condition. Their relative equilibria can be easily
characterized in the cases $|\Gamma_1|=|\Gamma_2| \equiv \Gamma$ .
Then we obtain that
$$ \| \vec c \|^2 + \| \vec w \|^2 = \Gamma \, \, \, \, \mbox{and} \, \, \, \,
\vec c \cdot \vec w = 0 .
$$
First we characterize the solutions such that   $\| \vec r_1(t) -
\vec r_2(t) \|$ is constant along the motion. Since
$$  \| \vec r_1 - \vec r_2 \|^2 = \| \vec c \|^2 \, \alpha_1^2 + \| \vec w \|^2 \, \alpha_2^2
  $$
 and since for$|\Gamma_1|=|\Gamma_2| $ either $\alpha_1 =0 $ or
 $\alpha_2 =0 $ it follows that
 $\| \vec c \|^2$ and $\| \vec w \|^2$ are constants for the relative equilibria.
 Therefore
 \begin{equation}
 \label{ref1} \vec c \cdot \dot {\vec c} = \epsilon \, \Delta_2 \, \vec c \cdot
( \vec w \times \hat z) = 0 .\end{equation}
 $\Delta_2 = 0$ implies that $w_z$ is a constant and therefore $\dot
 w_z = 0$. From the equations it follows that $\vec c \cdot ( \vec w \times \hat z) =
 0.$ This implies that the only way (\ref{ref1}) can be zero is if $ \vec c \cdot ( \vec w \times \hat z) =
 0.$ We arrive at the result that in a relative equilibria $\| \vec c \|^2$, $\| \vec w
 \|^2$ and $w_z$ are constants and that the vectors $\vec c \, , \,
 \vec w $ and $\hat z$ are always in the same plane. The vectors
 $\vec c$ and $\vec w$ will precess around the $z$-axis with
 constant angular velocity and keeping their relative orientation.
 It follows from the definition that $\vec r_1$ and $\vec r_2$
 will also precess around the $z$-axis with
 constant angular velocity and keeping their relative orientation.
 Therefore the only solutions for which $\| \vec r_1(t) -
\vec r_2(t) \|$ are constants are precisely the relative
equilibria.
\par

 We compute the angular velocity of the relative equilibria
 in the cases where $| \Gamma_1 | = | \Gamma_2 |$. \bigskip \par
 \noindent Case $ \Gamma_1 = \Gamma_2 $: If $\Gamma_1 =
 \Gamma_2 = \Gamma $ then $\alpha_1 = 0 $ and $\alpha_2 = \Gamma^{-1} \, .$
 Assuming that $c_x =  A \, \cos( \Omega \, t )$,
 $ c_y =  A \, \sin( \Omega \, t )$, $w_x = B \,  \cos( \Omega \, t ) $
 and $w_y = B \, \sin(\Omega \, t)$ we obtain that the angular velocity is given by
 $$ \Omega = \frac{\Gamma^2}{\| \vec w \|^2}\, \left( c_z - \frac{\sqrt{\| \vec c\|^2 - c_z^2}}{\mu} \right)
 + \frac{\epsilon}{2 \, \Gamma^2} \, \left( c_z + w_z \, \frac{\sqrt{\| \vec c\|^2 - c_z^2}}{\mu}
 \, \right) \, .$$
We observe that the perturbation factor in the angular velocity is
proportional to $\displaystyle \frac{1}{\Gamma^2}$ and therefore
the deviation from the corresponding spherical relative equilibria
should be bigger for smaller vorticities. We also observe that
when $\epsilon = 0 $ then, taking as before $\| \vec c \| = | c_z
| $, the above expression reduces to (\ref{gamma}). \bigskip
 \par
 \noindent Case $ \Gamma_1 = -\Gamma_2 $: If $\Gamma_1 =
 -\Gamma_2 = \Gamma $ then $\alpha_1 = \Gamma^{-1} $ and $\alpha_2 = 0 \, .$
 Proceeding as before we obtain
$$ \Omega = \frac{\Gamma^2}{\| \vec c \|^2}\, \left( c_z - \frac{\sqrt{\| \vec c\|^2 - c_z^2}}{\mu} \right)
 + \frac{\epsilon}{2 \, \Gamma^2} \, \left( c_z + w_z \, \frac{\sqrt{\| \vec c\|^2 - c_z^2}}{\mu}
 \, \right) \, .$$
The perturbation term is equal in both cases and again, as
$\epsilon \to 0 $ we will have $w \to \gamma \, .$ \par

\subsection{The 3-Vortex problem}
 The symmetry breaking when $\epsilon \ne 0$ is
expected to destroy the integrability of the three-vortices
spherical problem. In fact, when $\epsilon=0$ the center of
vorticity vector $ \vec c = \Gamma_1 \, \vec r_1 +
 \Gamma_2 \, \vec r_2 + \Gamma_3 \, \vec r_3 $ is conserved. The integrability
 follows by dimensionality. When
 $\epsilon \ne 0 $, only it's z-component $c_z =  \Gamma_1 \, z_1 +
 \Gamma_2 \, z_2 + \Gamma_3 \, z_3$ is conserved.
  That the integrability is destroyed
  is shown numerically in this section. We do
not pursue a complete numerical study of the problem. This will
addressed in a future work. \par
 To calculate the Poincar\'e sections for the three-vortex problem we
  will introduce an appropriate coordinate system.
First we introduce a cylindrical system of coordinates $(z_1 , z_2,
z_3 , \phi_1 , \phi_2 , \phi_3)$ with $\phi_i \in [ 0 , 2 \, \pi )$
and $z_i \in (-1,1)$ for $i=1,2$. The transformation is defined
trough
$$\begin{array}{l}
x_i = \sqrt{1-z_i^2} \, \cos(\phi_i) \, , \\
\\
y_i = \sqrt{1-z_i^2} \, \sin(\phi_i) \, , \\
\\
z_i = z_i \, . \end{array} $$ Since
$$w = \sum_{i=1}^3 \Gamma_i \, x_i \, dy_i \wedge dz_i + y_i \, dz_i \wedge dx_i
+ z_i \, dx_i \wedge dy_i = \sum_{i=1}^3 \Gamma_i \, d \phi_i
\wedge d z_i \, ,$$ this is a symplectic transformation and the
equations in cylindrical coordinates are given by
$$ \dot \phi_i = \frac{1}{\Gamma_i}\, \frac{\p H}{\p z_i} \, \, \, \mbox{ and }
 \dot z_i = - \frac{1}{\Gamma_i}\, \frac{\p  H}{\p \phi_i } \, \, \mbox{for} \, \, i=1,2 \, . $$
 To reduce the symmetry we define the coordinates
 $$ \begin{array}{ccl}

 q_1 & = &\Gamma_1 \, z_1 - \Gamma_3 \, z_3 \, ,
 \\
 \\
 q_2 & = & \Gamma_1 \, z_1 - \Gamma_2 \, z_2 \, ,
 \\
 \\
 q_3 & = & \Gamma_1 \, z_1 + \Gamma_2 \, z_2 + \Gamma_3 \, z_3 \, ,
 \\
 \\
 p_1 & = & \frac{\phi_1}{3} + \frac{\phi_2}{3}
 - 2 \,  \frac{\phi_3}{3} \,
 \\
 \\
p_2 & = & \frac{\phi_1}{3} -2 \,  \frac{\phi_2}{3} +
 \frac{\phi_3}{3} \,,
 \\
 \\
 p_3 & = &  \frac{\phi_1}{3} + \frac{\phi_2}{3} +
 \frac{\phi_3}{3} \, .
 \\
 \end{array} $$
 It follows that
 $$ dq_1 \wedge dp_1 + dq_2 \wedge dp_2 + dq_3 \wedge dp_3 = \Gamma_1 \, dz_1 \wedge d\phi_1 +
 \Gamma_2 \, dz_2 \wedge d\phi_2 + \Gamma_3 \, dz_3 \wedge d\phi_3 \, . $$

 Therefore the new system of coordinates forms a canonical system.
We observe that $q_3(t) = c_z$ and therefore $\dot q_3(t) =
\frac{\p H}{\p p_3}=0$. Then $p_3$ is a cyclic variable and the
Hamiltonian (\ref{hamN}) is given by
$$ H =H(q_1,q_2,p_1,p_2,c_z). $$ Since this Hamiltonian is now defined
in a four dimensional phase-space, the Poincar\`e maps can be
visualized. The equations are integrated numerically using the
Fehlberg-7(8) Runge-Kuta. The relative error in energy is kept
smaller then $10^{-10}$. A constant step size equals to $10^{-3}$
is used in all integrations. \par Figure 1 show the pictures of
three Poincar\`e maps. The onset of chaotic behavior is clear. It
appears (as expected) around the homoclinic solution present in
the integrable case $\epsilon = 0.00$. We observe that for
$\epsilon =0.05 $ and $\epsilon = 0.10$ the eccentricity of the
generating ellipses are $0.22$ and $0.33$ respectively.

 \end{document}